\documentclass[11pt,reqno]{amsart}
\usepackage{amssymb,amscd,amsbsy}
\usepackage{amssymb,amscd,amsbsy,mathrsfs}
\newcommand{\lb}{\linebreak}

\renewcommand{\a}{\alpha}

\newcommand{\z}{\zeta}

\renewcommand{\t}{\tau}

\newcommand{\f}{\varphi}
\renewcommand{\o}{\omega}

\renewcommand{\L}{\Lambda}

\newcommand{\h}{{\mathscr H}}

\newcommand{\K}{{\mathscr K}}

\newcommand{\T}{{\Bbb T}}

\newcommand{\dd}{{\Bbb D}}
\newcommand{\R}{{\Bbb R}}
\newcommand{\Z}{{\Bbb Z}}

\newcommand{\bS}{{\boldsymbol S}}

\newcommand{\rf}[1]{(\ref{#1})}

\newcommand{\df}{\stackrel{\mathrm{def}}{=}}

\newcommand{\supp}{\operatorname{supp}}

\newcommand{\const}{\operatorname{const}}

\newcommand{\eeq}{\end{equation}}
\newcommand{\beq}{\begin{equation}}
\newcommand{\bay}{\begin{eqnarray}}
\newcommand{\ba}{\begin{align*}}
\newcommand{\ea}{\end{align*}}
\newcommand{\ey}{\end{eqnarray}}
\newcommand{\bey}{\begin{eqnarray*}}
\newcommand{\eey}{\end{eqnarray*}}

\newcommand{\be}{\infty}

\newcommand{\bl}{\blacksquare}

\newcommand{\Pf}{{\bf Proof. }}

\newcommand{\ov}{\overline}

\newtheorem{thm}{\hspace{\parindent}Theorem}[section]

\newtheorem{cor}[thm]{\hspace{\parindent}Corollary}
\newtheorem{lem}[thm]{\hspace{\parindent}Lemma}

\theoremstyle{remark}

\newtheorem*{rem*}{Remark}

\newcommand\CA{{\rm C}_{\rm A}}

\begin{document}

\newcommand{\vse}{\vspace{.2in}}
\numberwithin{equation}{section}

\title{Functions of commuting contractions under perturbation}
\author{V.V. Peller}

\begin{abstract}
The purpose of the paper is to obtain estimates for differences of functions of two pairs of commuting contractions on Hilbert space. In particular, Lipschitz type estimates, H\"older type estimates, Schatten--von Neumann estimates are obtained.
The results generalize earlier known results for functions of self-adjoint operators, normal operators, contractions and dissipative operators.
\end{abstract}

\maketitle


\setcounter{section}{0}
\section{\bf Introduction}
\label{In}

\

We are going to obtain sharp estimates for the differences
$f(T_1,R_1)-f(T_2,R_2)$ for functions $f$ analytic on the bidisk and 
for pairs $(T_1,R_1)$ and $(T_2,R_2)$ of commuting contractions on Hilbert space.
We obtain analogs of earlier results for functions of self-adjoint operators, normal operators, commuting $n$-tuples of self-adjoint operators, contractions and dissipative operators, see \cite{Pe2}, \cite{Pe3}, \cite{Pe+}, \cite{AP1}, \cite{AP2}, \cite{Pe*}, \cite{AP3}, \cite{APPS}, \cite{NP}, \cite{AP5}. We also refer the reader to recent surveys \cite{Pe6} and \cite{AP4}.

Let $(T,R)$ be a pair of commuting contractions (i.e., operators of norm at most one) on a Hilbert space $\h$.
By Ando's theorem \cite{An} (see also \cite{SNF}, Ch. I, \S\:6), $T$ and $R$ have commuting unitary dilations, i.e., there exist commuting unitary operators $U$ and $V$ on a Hilbert space $\K$, $\K\supset\h$, such that
$$
T^jR^k=PU^jV^k\big|\h,\quad j,~k\in\Z_+,
$$
where $P$ is the orthogonal projection onto $\h$.
Ando's theorem allows us construct a functional calculus on the bidisk algebra 
$\CA(\dd^2)$, which consists of functions of two variables analytic in the bidisk $\dd^2$ and continuous in its closure. Indeed, for $f\in\CA(\dd^2)$, we put
$$
f(T,R)\df Pf(U,V)\big|\h.
$$
This functional calculus is linear and multiplicative, and  the following analog of von Neumann's inequality holds:
$$
\|f(T,R)\|\le\max\{|f(\z_1,\z_2)|:~|\z_1|\le1,~|\z_2|\le1\},\quad f\in\CA(\dd^2).
$$

We are going to obtain the following Lipschitz type estimate:
\bay
\label{oLio2szh}
\|f(T_1,R_1)-f(T_2,R_2)\|\le\const\max\big\{\|T_1-T_2\|,\|R_1-R_2\|\big\}
\ey
for pairs $(T_1,R_1)$ and $(T_2,R_2)$ of commuting contractions; this inequality
holds for analytic functions $f$ in the Besov space $B_{\be,1}^1(\T^2)$. We are also going to obtain H\"older type estimates, Schatten--von Neumann estimates and related results.

Note that Farforovskaya showed in \cite{F} that Lipschitz functions on the real line $\R$ are not necessarily {\it operator Lipschitz}, i.e., the condition 
$|f(x)-f(y)|\le\const|x-y|$ does not imply that 
$$
\|f(A)-f(B)\|\le\const\|A-B\|
$$
for arbitrary bounded self-adjoint operators $A$ and $B$. The same is true for functions of unitary operators, see e.g., \cite{AP4}. 

Later, it was shown in \cite{Pe2} and \cite{Pe3} that if $f$ belongs to the (homogeneous) Besov space $B_{\be,1}^1(\R)$, then $f$ is operator Lipschitz.
We obtain in \S\:\ref{Lityes} an analog of this result for functions of two commuting contractions. 

It was also shown in \cite{Pe2} that if a function $f$ on the unit circle is operator Lipschitz, then it must belong to the Besov space $B_{1,1}^1(\T)$ (this was deduced from the trace class criterion for Hankel operators, see \cite{Pe1} and \cite{Pe4}).
This immediately implies that inequality \rf{oLio2szh} does not hold for all analytic continuously differentiable functions on $\T^2$.

It was discovered in \cite{AP1} that the situation with H\"older type estimates is quite different from Lipschitz type estimates. It was shown in \cite{AP1} that H\"older functions  of order $\a$, $0<\a<1$, on $\R$ and on $\T$ are necessarily operator H\"older of order $\a$. In \S\:\ref{Hotyespromone} we establish H\"older type operator estimates  for functions of pairs of commuting contractions.

In \S\:\ref{otsenkiSchsvNe} we obtain Schatten--von Neumann estimates for the differences $f(T_1,R_1)-f(T_2,R_2)$, where $(T_1,R_1)$ and $(T_2,R_2)$
are pairs of commuting contractions. Finally, \S\:\ref{komu tator}
is devoted to estimates of (quasi)commutators $f(T_1,R_1)Q-Qf(T_2,R_2)$.

To obtain such estimates in the case of functions of self-adjoint operators, normal operators, contractions, dissipative operators , techniques of double operator integrals were used (see \cite{Pe2}, \cite{Pe3}, \cite{Pe+}, \cite{AP1}, \cite{AP2}, \cite{Pe*}, \cite{AP3}, \cite{APPS}, \cite{NP}, \cite{AP5}). Double operator integrals play an important role in perturbation theory. They appeared in the paper \cite{DK} by Daletskii and Krein. Later Birman and Solomyak developed a beautiful theory of double operator integrals in \cite{BS1}, \cite{BS2} and \cite{BS3}.

To prove the results of this paper, we could also use double operator integrals. However, I have decided to make the presentation elementary and avoid double operator integrals in the proofs by replacing double operator integrals with elementary identities. Note that these identities are inspired by double operator integrals. Nevertheless, to understand the proofs, no knowledge of double operator integrals is assumed. 

I would like to express my gratitude to A.B. Alexandrov and M.M. Malamud for helpful discussions.

\

\section{\bf Besov classes of periodic functions}
\label{Bes}

\

In this section we give a brief introduction to Besov spaces on the torus.

To define Besov spaces on the torus $\T^d$, we consider an infinitely differentiable function $w$ on $\R$ such
that
$$
w\ge0,\quad\supp w\subset\left[\frac12,2\right],\quad\mbox{and} \quad w(s)=1-w\left(\frac s2\right)\quad\mbox{for}\quad s\in[1,2].
$$
Let $W_n$, $n\ge0$, be the trigonometric polynomials defined by
$$
W_n(\z)\df\sum_{j\in\Z^d}w\left(\frac{|j|}{2^n}\right)\z^j,\quad n\ge1,
\quad W_0(\z)\df\sum_{\{j:|j|\le1\}}\z^j,
$$
where 
$$
\z=(\z_1,\cdots,\z_d)\in\T^d,\quad j=(j_1,\cdots,j_d),\quad\mbox{and}\quad
|j|=\big(|j_1|^2+\cdots+|j_d|^2\big)^{1/2}.
$$
For a distribution $f$ on $\T^d$ we put
\bay
\label{fnWn}
f_n=f*W_n,\quad n\ge0.
\ey
It is easy to see that
\bay
\label{fSigmafn}
f=\sum_{n\ge0}f_n;
\ey
the series converges in the sense of distributions.
We say that $f$ belongs the {\it Besov class} $B_{p,q}^s(\T^d)$, $s>0$, 
$1\le p,\,q\le\be$, if
\bay
\label{Bperf}
\big\{2^{ns}\|f_n\|_{L^p}\big\}_{n\ge0}\in\ell^q.
\ey
The Besov class $B_{\be,\be}^\a(\T^d)$, $\a>0$, coincides with the class $\L_\a(\T^d)$
of H\"older--Zygmund functions of order $\a$. If $0<\a<1$, then $\L_\a(\T^d)$
is the class of functions  $f$ on $\T^d$ satisfying the inequality
$$
|f(\z_1,\cdots,\z_d)-f(\t_1,\cdots,\t_d)|\le
\const\Big(\max_{1\le j\le d}|\z_j-\t_j|\Big)^\a.
$$

The analytic subspace $B_{p,q}^s(\T^d)_+$ of $B_{p,q}^s(\T^d)$ consists of functions $f$ in $B_{p,q}^s(\T^d)$ for which the Fourier coefficients 
$\widehat f(j_1,\cdots,j_d)$ satisfy the equalities:
\bay
\label{anafunapo}
\widehat f(j_1,\cdots,j_d)=0\quad\mbox{whenever}\quad
\min_{1\le k\le d}j_k<0.
\ey
We use the notation $\L_\a(\T^d)_+$ for $B_{\be,\be}^\a(\T^d)_+$, $\a>0$.

We refer the reader to \cite{Pee} for more detailed information about Besov spaces.

\

\section{\bf Certain operator inequalities}
\label{nerva}

\

We are going to use in this section Lemma 3.2 of \cite{AP4}. Suppose that $\{A_j\}_{j\ge0}$ is a sequence of bounded linear operators on Hilbert space such that
\bay
\label{nrvadlyaAj}
\left\|\sum_{j\ge0}A_j^*A_j\right\|\le1\quad\mbox{and}\quad
\left\|\sum_{j\ge0}A_jA_j^*\right\|\le1.
\ey
Let $Q$ be a bounded linear operator. Consider the row ${\rm R}_{\{A_j\}}(Q)$ and
the column ${\rm C}_{\{A_j\}}(Q)$ defined by
$$
{\rm R}_{\{A_j\}}(Q)\df\big(A_0Q\:A_1Q\:A_2Q\:\cdots\big)
$$
and
$$
{\rm C}_{\{A_j\}}(Q)\df\left(\begin{matrix}QA_0\\QA_1\\QA_2\\\vdots\end{matrix}\right).
$$
Then for $p\in[2,\be]$,
\bay
\label{RCAjSp}
\big\|{\rm R}_{\{A_j\}}(Q)\big\|_{\bS_p}\le\|Q\|_{\bS_p}\quad\mbox{and}\quad
\big\|{\rm C}_{\{A_j\}}(Q)\big\|_{\bS_p}\le\|Q\|_{\bS_p}
\ey
whenever $Q\in\bS_p$, see Lemma 3.2 of \cite{AP4}.

Note that in the case $p=\be$ we assume that $Q$ is an arbitrary bounded linear operator and by $\|\cdot\|_\be$ we mean operator norm. 

\begin{thm}
\label{Aj1QAj2}
Let $Q$ be a bounded linear operator and let $\big\{A_j^{(1)}\big\}_{j\ge1}$
and $\big\{A_j^{(2)}\big\}_{j\ge1}$ be sequences of bounded linear operators satisfying {\em\rf{nrvadlyaAj}}. Then the series
\bay
\sum_{j\ge1}A_j^{(1)}QA_j^{(2)},
\ey
converges in the weak operator topology and the sum satisfies the estimate
$$
\left\|\sum_{j\ge0}A_j^{(1)}QA_j^{(2)}\right\|_{\bS_p}\le\|Q\|_{\bS_p},
\quad1\le p\le\be,
$$
whenever $Q\in\bS_p$.
\end{thm}

\Pf We can represent $Q$ as $Q=Q_1Q_2$, where $Q_1,\,Q_2\in\bS_{2p}$ and 
$\|Q_1\|_{\bS_2p}=\|Q_2\|_{\bS_2p}=\|Q\|_{S_p}^{1/2}$. Clearly,
$$
\sum_{j\ge1}A_j^{(1)}QA_j^{(2)}=
{\rm R}_{\big\{A^{(1)}_j\big\}}(Q_1){\rm C}_{\big\{A^{(2)}_j\big\}}(Q_2), 
$$
and so by \rf{RCAjSp},
\begin{align*}
\left\|\sum_{j\ge1}A_j^{(1)}QA_j^{(2)}\right\|_{\bS_p}&\le
\Big\|{\rm R}_{\big\{A^{(1)}_j\big\}}(Q_1)\Big\|_{\bS_2p}
\Big\|{\rm C}_{\big\{A^{(2)}_j\big\}}(Q_2)\Big\|_{\bS_2p}\\[.2cm]
&\le\|Q_1\|_{\bS_{2p}}\|Q_2\|_{\bS_{2p}}=\|Q\|_{\bS_p},\quad Q\in\bS_p.
\quad\bl
\end{align*}

\begin{lem}
\label{AjfjTR}
Let $T$ and $R$ be commuting contractions on a Hilbert space $\h$.
Suppose that $f_1,f_2,\cdots,f_n$ are analytic polynomials of two variables
such that 
$$
\sum_{j=1}^n|f_j(\z,\t)|^2\le1,\quad |\z|,~|\t|\le1.
$$
Let $A_j=f_j(T,R)$, $0\le j\le n$. Then the operators $A_j$ satisfy
{\em\rf{nrvadlyaAj}}.
\end{lem}

\Pf Let $U$ and $V$ be commuting unitary dilations of $T$ and $R$ on a Hilbert space 
$\K$, $\K\supset\h$, and let $P$ be the orthogonal projection onto $\h$. Then
$$
A_jx=f_j(T,R)x=Pf_j(U,V)x,\quad x\in\h,
$$
and
$$
A^*_jx=(f_j(T,R))^*x=P\ov f_j(U,V)x,\quad x\in\h.
$$
Suppose that $E_{U,V}$ is the joint spectral measure of $U$ and $V$ on $\T^2$.
Let $x\in\h$. We have
\begin{align*}
\sum_{j=1}^n(A_j^*A_jx,x)&=\sum_{j=1}^n(A_jx,A_jx)=
\sum_{j=1}^n\big(Pf_j(U,V)x,Pf_j(U,V)x\big)\\[.2cm]
&\le\sum_{j=1}^n\big(f_j(U,V)x,f_j(U,V)x\big)=
\left(\left(\sum_{j=1}^n|f_j|^2(U,V)\right)x,x\right)\\[.2cm]
&=\int_{\T\times\T}\sum_{j=1}^n|f_j|^2\,(dE_{U,V}x,x)
\le\|x\|^2.
\end{align*}
On the other hand,
\begin{align*}
\sum_{j=1}^n(A_jA^*_jx,x)&=\sum_{j=1}^n(A^*_jx,A^*_jx)=
\sum_{j=1}^n\big(P\ov f_j(U,V)x,P\ov f_j(U,V)x\big)\\[.2cm]
&\le\sum_{j=1}^n\big(\ov f_j(U,V)x,\ov f_j(U,V)x\big)=
\left(\left(\sum_{j=1}^n|f_j|^2(U,V)\right)x,x\right)
\le\|x\|^2.
\end{align*}
This completes the proof. $\bl$

\begin{cor}
\label{rzberyo}
Suppose that $\f_1,\f_2,\cdots,\f_n$ and $\psi_1,\psi_2,\cdots,\psi_n$
are analytic polynomials of two variables. Put
$$
M_1\df\max_{(\z,\t)\in\dd^2}\sum_{j=1}^n|\f_j(\z,\t)|^2
\quad\mbox{and}\quad
M_2\df\max_{(\z,\t)\in\dd^2}\sum_{j=1}^n|\psi_j(\z,\t)|^2.
$$
Let $(T_1,R_1)$ and $(T_2,R_2)$ be pairs of commuting contractions and let
$Q$ be a bounded linear operator on Hilbert space. Then
$$
\left\|\sum_{j=1}^n\f_j(T_1,R_1)Q\psi_j(T_2,R_2)\right\|_{\bS_p}\le
(M_1M_2)^{1/2}\|Q\|_{\bS_p},\quad1\le p\le\be,
$$
whenever $Q\in\bS_p$.
\end{cor}

\Pf The result follows immediately from Theorem \ref{Aj1QAj2} and from Lemma
\ref{AjfjTR}. $\bl$

\

\section{\bf Lipschitz type estimates}
\label{Lityes}

\

We start this sections with an algebraic formula for 
$f(T_1,R_1)-f(T_2,R_2)$ for an analytic polynomial $f$ in two variables
and for pairs $(T_1,R_1)$ and $(T_2,R_2)$ of commuting contractions.
This identity plays a crucial role in what follows.

We define the linear operators $S_1^*$ and $S_2^*$ on the set of analytic functions on $\dd^2$.
For $f(z_1,z_2)=\sum_{k,m\ge0}a_{k,m}z_1^kz_2^m$, we put
$$
(S_1^*f)(z_1.z_2)\df\sum_{k,\,m\ge0}a_{k+1,m}z_1^kz_2^m\quad\mbox{and}\quad
(S_2^*f)(z_1.z_2)\df\sum_{k,\,m\ge0}a_{k,m+1}z_1^kz_2^m.
$$
The algebraic identity in the following theorem involves powers of the operators $S_1^*$ and $S_2^*$.

\begin{thm}
\label{tozddlyapol}
Let $f$ be an analytic polynomial of two variables of degree at most $n$ in each variable. Suppose that $(T_1,R_1)$ and $(T_2,R_2)$ are pairs of commuting contractions. Then the following identity holds:
\begin{align}
\label{predraz}
f(T_1,R_1)-f(T_2,R_2)&=\sum_{j=1}^{n}\big(\big((S_2^*)^jf\big)(T_1,R_1)\big)(R_1-R_2)R_2^{j-1}
\nonumber\\[.2cm]
&+\sum_{j=1}^{n}T_1^{j-1}(T_1-T_2)\big((S_1^*)^jf\big)(T_2,R_2).
\end{align}
\end{thm}

\Pf Clearly, it suffices to establish \rf{predraz} for monomials $f$, i.e., 
polynomials of the form $f(z_1,z_2)=z_1^lz_2^m$ with $l,\,m\le n$. We have
\begin{multline*}
\sum_{j=1}^{n}\big((S_2^*)^jf\big)(T_1,R_1)(R_1-R_2)R_2^{j-1}=
\sum_{j=1}^mT_1^lR_1^{m-j}(R_1-R_2)R_2^{j-1}\\[.2cm]
=\big(T_1^lR_1^m-T_1^lR_1^{m-1}R_2\big)
+\big(T_1^lR_1^{m-1}R_2-T_1^lR_1^{m-2}R_2^2\big)+\cdots\\[.2cm]
+\big(T_1^lR_1R_2^{m-1}-T_1^lR_2^m\big)=T_1^lR_1^m-T_1^lR_2^m.
\end{multline*}
Similarly,
$$
\sum_{j=1}^{n}T_1^{j-1}(T_1-T_2)\big((S_1^*)^jf\big)(T_2,R_2)=
\sum_{j=1}^lT_1^{j-1}(T_1-T_2)T_2^{l-j}R_2^m=-T_2^lR_2^m+T_1^lR_2^m.
$$
Thus, the right-hand side of \rf{predraz} is equal to
$$
T_1^lR_1^m-T_2^lR_2^m=f(T_1,R_1)-f(T_2,R_2).\quad\bl
$$

The following result is a Bernstein type inequality for functions of
commuting contractions.

%

\begin{thm}
\label{Lipoden}
Let $f$ be an analytic polynomial of two variables of degree at most $n$ in each variable. Then
\bay
\label{nervotB}
\|f(T_1,R_1)-f(T_2,R_2)\|\le\const n\|f\|_{L^\be(\T^2)}\max\big\{\|T_1-T_2\|,\|R_1-R_2\|\big\}.
\ey
\end{thm}

We need the following lemma.

\begin{lem}
\label{S*fj}
Let $f$ be a function on the unit disk $\dd$ of class $H^\be$. Then
$$
\left(\sum_{j=1}^n\big|\big((S^*)^jf\big)(\z)\big|^2\right)^{1/2}\le\const n^{1/2}\|f\|_{H^\be},\quad\mbox{a.e. on}\quad\T.
$$
\end{lem}

Here $S^*$ is the operator on the space of functions analytic in $\dd$ defined by 
$$
S^*\left(\sum_{j\ge0}a_jz^j\right)=\sum_{j\ge0}a_{j+1}z^j.
$$

{\bf Proof of Lemma \ref{S*fj}.} Let $\varSigma_jf$ be the $j$th partial sum of
the Taylor series of $f$:
$$
(\varSigma_jf)(z)=\sum_{k=0}^j\widehat f(k)z^k.
$$
Clearly, $f=\varSigma_{j-1}f+z^j(S^*)^jf$, and so it suffices to show that
$$
\left(\sum_{j=1}^n\big|\big(\varSigma_{j-1}f\big)(\z)\big|^2\right)^{1/2}\le\const n^{1/2}\|f\|_{H^\be},\quad|\z|=1.
$$
Clearly, it suffices to prove the last inequality for $\z=1$. Consider the function $F$ analytic in $\dd$ defined by $F=(1-z)^{-1}f(z)$. Clearly,
$\widehat F(j)=(\varSigma_jf)(1)$. Thus, it suffices to show that
\bay
\label{L2chsu}
\left(\sum_{k=0}^n|\widehat F(k)|^2\right)^{1/2}\le\const n^{1/2}\|f\|_{H^\be}.
\ey
Let $F_r(z)\df F(rz)$, $0<r<1$. It is easy to see that
$$
\|F_r\|_{L^2(\T)}\le\left\|\frac1{1-rz}\right\|_{L^2(\T)}\|f\|_{H^\be}
=(1-r^2)^{-1/2}\|f\|_{H^\be}.
$$
Clearly, we can assume that $n\ge2$. Put $r=1-1/n$. We have
$$
\sum_{k=0}^n|\widehat F(k)|^2\le\frac1{r^{2n}}\sum_{k=0}^nr^{2k}|\widehat F(k)|^2
\le\frac1{r^{2n}}\|F_r\|^2_{L^2\T)}\le\frac{(1-r^2)^{-1}}{r^{2n}}\|f\|^2_{H^\be}
\le\const n\|f\|^2_{H^\be}
$$
which proves \rf{L2chsu}. $\bl$

\medskip

{\bf Proof of Theorem \ref{Lipoden}.} We are going to use formula \rf{predraz}
end estimate the norms of each term of the right-hand side. Put 
$\f_j\df(S_2^*)f$ and $\psi_j(z_1,z_2)\df z_2^{j-1}$, $1\le j\le n$.
Then 
$$
\sum_{j=1}^{n}\big(\big((S_2^*)^jf\big)(T_1,R_1)\big)(R_1-R_2)R_2^{j-1}
=\sum_{j=1}^{n}\f_j(T_1,R_1)(R_1-R_2)\psi_j(T_2,R_2).
$$
It follows now from Corollary \ref{rzberyo}, Lemma \ref{S*fj} and the trivial equality
$$
\sum_{j=1}^n|\psi_j(\z_1,\z_2)|^2=n,\quad|\z_1|=|\z_2|=1,
$$
that
$$
\left\|\sum_{j=1}^{n}\big(\big((S_2^*)^jf\big)(T_1,R_1)\big)(R_1-R_2)R_2^{j-1}
\right\|\le\const\cdot n\|f\|_{L^\be(\T^2)}\|R_1-R_2\|.
$$
In the same way it can be shown that
$$
\left\|\sum_{j=1}^{n}T_1^{j-1}(T_1-T_2)\big((S_1^*)^jf\big)(T_2,R_2)\right\|
\le\const\cdot n\|f\|_{L^\be(\T^2)}\|T_1-T_2\|.\quad\bl
$$

We are ready to proceed to the main result of this section.

\begin{thm}
\label{LiotsfBes}
Let $(T_1,R_1)$ and $(T_2,R_2)$ be pairs of commuting contractions. Suppose that $f$ is a function in two variables of class $B_{\be,1}^1(\T^2)_+$. Then
$$
\|f(T_1,R_1)-f(T_2,R_2)\|
\le\const\|f\|_{B_{\be,1}^1(\T^2)}\max\big\{\|T_1-T_2\|,\|R_1-R_2\|\big\}.
$$
\end{thm}

\Pf Let $f_n$ be the analytic polynomial defined by \rf{fnWn}. Clearly, $f_n$ has degree at most $2^{n+1}-1$ in each variable. It follows from \rf{fSigmafn} and from
Theorem \ref{Lipoden} that
\begin{align*}
\|f(T_1,R_1)-f(T_2,R_2)\|&\le\sum_{n\ge0}\|f_n(T_1,R_1)-f_n(T_2,R_2)\|\\[.2cm]
&\le\const\sum_{n\ge0}n\|f_n\|_{L^\be}\max\big\{\|T_1-T_2\|,\|R_1-R_2\|\big\}\\[.2cm]
&\le\const\|f\|_{B_{\be,1}^1}\max\big\{\|T_1-T_2\|,\|R_1-R_2\|\big\}
\end{align*}
(see the definition of the Besov class $B_{\be,1}^1(\T^2)_+$ given in \S\:\ref{Bes}).
$\bl$

\

\section{\bf H\"older type estimates and arbitrary moduli of continuity}
\label{Hotyespromone}

\

In \cite{AP1} Bernstein type inequalities for functions of self-adjoint and unitary operators were used to obtain H\"older type estimates for functions of self-adjoint and unitary operators. The same method allows us to deduce from Bernstein type inequality \rf{nervotB} the following H\"older type estimates for functions 
of pairs of commuting contractions.

\begin{thm}
Let $0<\a<1$ and let $f$ be an analytic function of class $\L_\a(\T^2)_+$.
Then 
$$
\|f(T_1,R_1)-f(T_2,R_2)\|
\le\const\|f\|_{\L_\a(\T^2)}\max\big\{\|T_1-T_2\|^\a,\|R_1-R_2\|^\a\big\}
$$
for arbitrary pairs $(T_1,R_1)$ and $(T_2,R_2)$ of commuting contractions.
\end{thm}

The theorem can be deduced from Theorem \ref{Lipoden} in the same way as it was done in the proof of Theorem 5.1 of \cite{AP1} in the case of functions of unitary operators.

Consider now an arbitrary modulus of continuity $\o$, i.e., a continuous nondecreasing function 
$\o$ on $[0,\be)$ such $\o(0)=0$ and $\o(s+t)\le\o(s)+\o(t)$, $s,\,t\ge0$.
We define the functions $\o_*$ by
$$
\o_*(s)=s\int_s^\be\frac{\o(t)}{t^2}\,dt=\int_1^\be\frac{\o(st)}{t^2}\,dt
,\quad s\ge0.
$$
It is easy to verify that if $\o_*(s)<\be$ for some $s>0$, then $\o_*(s)<\be$ for every $s>0$ and $\o_*$ is a modulus of continuity, see \cite{AP1}.

We use the notation $\L_\o(\T^2)_+$ for the class of analytic functions $f$
on $\T^2$ (i.e., functions $f$ satisfying \rf{anafunapo}) such that
$$
|f(\z_1,\t_1)-f(\z_2,\t_2)|\le\const\o\big(\max\{|\z_1-\z_2|,|\t_1-\t_2|\}\big),
\qquad(\z_1,\t_1),~(\z_2,\t_2)\in\T^2.
$$

\begin{thm}
Let $\o$ be a modulus of continuity and let $f\in\L_\o(\T^2)_+$. Then
$$
\|f(T_1,R_1)-f(T_2,R_2)\|
\le\const\|f\|_{\L_\o(\T^2)}
\o_*\big(\max\{\|T_1-T_2\|,\|R_1-R_2\|\}\big)
$$
for arbitrary pairs $(T_1,R_1)$ and $(T_2,R_2)$ of commuting contractions.
\end{thm}

Again, the theorem can be deduced from Theorem \ref{Lipoden} in the same way as it was done in the proof of Theorem 7.1 of \cite{AP1} in the case of functions of unitary operators.

\

\section{\bf Schatten--von Neumann estimates}
\label{otsenkiSchsvNe}

\

We obtain in this sections Schatten--von Neumann estimates for the differences
\lb$f(T_1,R_1)-f(T_2,R_2)$, where $(T_1,R_1)$ and $(T_2,R_2)$ 
are pairs of commuting contractions. The following theorem can be deduced from 
Corollary \ref{rzberyo} and Theorem \ref{tozddlyapol} in the same way as in the proof of Theorem \ref{LiotsfBes}.

\begin{thm}
\label{SpLiotsfBes}
Let $1\le p<1$ and
let $(T_1,R_1)$ and $(T_2,R_2)$ be pairs of commuting contractions
such that $T_1-T_2\in\bS_p$ and $R_1-R_2\in\bS_p$. Suppose that $f$ is a function of class $B_{\be,1}^1(\T^2)_+$. Then $f(T_1,R_1)-f(T_2,R_2)\in\bS_p$ and
$$
\|f(T_1,R_1)-f(T_2,R_2)\|_{\bS_p}\le\const
\|f\|_{B_{\be,1}^1(\T^2)}\max\big\{\|T_1-T_2\|_{\bS_p},\|R_1-R_2\|_{\bS_p}\big\}.
$$
\end{thm}

We proceed now to Schatten--von Neumann estimates for $f(T_1,R_1)-f(T_2,R_2)$ in the case when $f$ belongs to the H\"older class $\L_\a(\T^2)_+$, $0<\a<1$.

\begin{thm}
\label{SpHootsfBes}
Let $1<p<1$ and
let $(T_1,R_1)$ and $(T_2,R_2)$ be pairs of commuting contractions
such that $T_1-T_2\in\bS_p$ and $R_1-R_2\in\bS_p$. Suppose that $f$ is a function  of class $\L_\a(\T^2)_+$, $0<\a<1$. Then $f(T_1,R_1)-f(T_2,R_2)\in\bS_{p/\a}$ and
\bay
\label{Honerp>1}
\|f(T_1,R_1)-f(T_2,R_2)\|_{\bS_{p/\a}}\le\const
\|f\|_{\L_\a(\T^2)}\max\big\{\|T_1-T_2\|^\a_{\bS_p},\|R_1-R_2\|^\a_{\bS_p}\big\}.
\ey
\end{thm}

The theorem can be deduced from Theorem \ref{SpLiotsfBes} practically in the same way as it was done in the proof of Theorem 6.7 of \cite{AP2}.

Note that the conclusion of Theorem \ref{SpHootsfBes} does not hold for $p=1$. To get inequality \rf{Honerp>1} with $p=1$, we can impose the condition $f\in B_{\be,1}^\a(\T^2)_+$,
see Theorem 6.3 of \cite{AP2} which is a similar result for functions of unitary operators. Also, let me mention that under the assumptions $f\in\L_\a(\T^2)_+$
and $p=1$, one can conclude that the singular values of $f(T_1,R_1)-f(T_2,R_2)$
satisfy the estimate
$$
s_j(f(T_1,R_1)-f(T_2,R_2))\le\const(1+j)^{-\a/p},\quad j\ge0,
$$
Again, this was proved in Theorem 6.1 of \cite{AP2} in the case of functions of unitary operators. In our situation the proofs are similar.

\

\section{\bf Commutator estimates}
\label{komu tator}

\

In this section we obtain estimates for (quasi)commutators of the form
$f(T_1,R_1)Q-Qf(T_2,R_2)$, where $(T_1,R_1)$ and $(T_2,R_2)$ are  pairs of commuting contractions and $Q$ is a bounded linear operator. In the special case $T_1=T_2$ and $R_1=R_2$, we arrive at commutator estimates for $f(T,R)Q-Qf(T,R)$, where $(T,R)$ is a pair of commuting contractions.

Let us first establish an analog of formula \rf{predraz} for quasicommutators.

\begin{thm}
Let $f$ be an analytic polynomial of two variables of degree at most $n$ in each variable. Suppose that $(T_1,R_1)$ and $(T_2,R_2)$ are pairs of commuting contractions and $Q$ is a bounded linear operator. Then the following identity holds:
\begin{align}
\label{predrazkom}
f(T_1,R_1)Q-Qf(T_2,R_2)&=
\sum_{j=1}^{n}\big(\big((S_2^*)^jf\big)(T_1,R_1)\big)(R_1Q-QR_2)R_2^{j-1}
\nonumber\\[.2cm]
&+\sum_{j=1}^{n}T_1^{j-1}(T_1Q-QT_2)\big((S_1^*)^jf\big)(T_2,R_2).
\end{align}
\end{thm}

The proof of formula \rf{predrazkom} is similar to the proof of formula \rf{predraz}.

We can deduce now from formula \rf{predrazkom} the following analog of Theorem 
\ref{LiotsfBes}.

\begin{thm}
\label{LiotsfBesko}
Let $(T_1,R_1)$ and $(T_2,R_2)$ be pairs of commuting contractions
and let $Q$ be a bounded linear operator. Suppose that $f$ is a function in two variables of class $B_{\be,1}^1(\T^2)_+$. Then
$$
\|f(T_1,R_1)Q-Qf(T_2,R_2)\|
\le\const\|f\|_{B_{\be,1}^1(\T^2)}\max\big\{\|T_1Q-QT_2\|,\|R_1Q-QR_2\|\big\}.
$$
\end{thm}

The proof of Theorem \ref{LiotsfBesko} is practically the same as the proof of Theorem \ref{LiotsfBes}.

We can also obtain analogs of the results of \S\:\ref{Hotyespromone} and 
\S\:\ref{otsenkiSchsvNe} of this paper for quasicommutators.
Let us state, for example, the following result, which can be proved in the same way as Theorem \ref{SpHootsfBes}.

\begin{thm}
Suppose that $1<p<1$ and $0<\a<1$.
Let $(T_1,R_1)$ and $(T_2,R_2)$ be pairs of commuting contractions
and let $Q$ be a bounded linear operator
such that $T_1Q-QT_2\in\bS_p$ and $R_1Q-QR_2\in\bS_p$. Suppose that $f$ is a function  of class $\L_\a(\T^2)_+$. Then $f(T_1,R_1)Q-Qf(T_2,R_2)\in\bS_{p/\a}$ and
\begin{multline*}
\label{Honerp>1}
\|f(T_1,R_1)Q-Qf(T_2,R_2)\|_{\bS_{p/\a}}\\[.2cm]
\le\const
\|f\|_{\L_\a(\T^2)}\max\big\{\|T_1Q-QT_2\|^\a_{\bS_p},\|R_1Q-QR_2\|^\a_{\bS_p}\big\}
\|Q\|^{1-\a}.
\end{multline*}
\end{thm}

\

\

\noindent
Department of Mathematics\\
Michigan State University\\
East Lansing Michigan 48824\\

\end{document}